\documentclass[a4paper,12pt]{article}

\usepackage{color}
\usepackage{pstricks}
\usepackage{pst-node}
\usepackage{makeidx}
\usepackage{latexsym}
\usepackage{amsmath}
\usepackage{amssymb}
\usepackage{amsthm}
\usepackage{amscd}
\usepackage{float}
\usepackage{psfrag}

\newcommand{\flata}{\ensuremath{\alpha}}

\newcommand{\canonicala}{\ensuremath{i}}
\newcommand{\canonicalb}{\ensuremath{j}}
\newcommand{\Gr}{\ensuremath{\mathrm{Gr}}}

\newcommand{\frakS}{\mathfrak{S}}
\newcommand{\rt}{\mathrm{right}}
\newcommand{\lt}{\mathrm{left}}

\newcommand{\sgn}{\ensuremath{\mathop{\mathrm{sgn}}}\nolimits}

\newcommand{\combination}[2]{\ensuremath{
 \begin{pmatrix} #1 \\ #2 \end{pmatrix}}}

\newcommand{\formal}{\mathrm{formal}}
\theoremstyle{plain}
\newtheorem{theorem}{Theorem}[section]

\newtheorem{corollary}[theorem]{Corollary}
\newtheorem{lemma}[theorem]{Lemma}

\newtheorem{conjecture}[theorem]{Conjecture}
\theoremstyle{definition}
\newtheorem{definition}[theorem]{Definition}

\newcommand{\suchthat}{\; | \;}
\renewcommand{\Re}{\ensuremath{\mathfrak{Re}}}

\newcommand{\diag}{\mathrm{diag}}
\newcommand{\Ext}{\mathop{\mathrm{Ext}}\nolimits}
\newcommand{\Hom}{\mathop{\mathrm{Hom}}\nolimits}

\newcommand{\GL}{GL}
\newcommand{\bCx}{\bC^{\times}}

\newcommand{\coh}{\mathrm{coh}}
\newcommand{\DbX}{D^b \mathrm{coh}(X)}

\renewcommand{\det}{\mathop{\mathrm{det}}\nolimits}

\newcommand{\bC}{\ensuremath{\mathbb{C}}}

\newcommand{\bP}{\ensuremath{\mathbb{P}}}

\newcommand{\bZ}{\ensuremath{\mathbb{Z}}}

\newcommand{\scE}{\ensuremath{\mathcal{E}}}

\newcommand{\scO}{\ensuremath{\mathcal{O}}}

\title{Stokes Matrices for the Quantum Cohomologies of Grassmannians}
\author{Kazushi Ueda}
\date{}
\begin{document}

\maketitle

\begin{abstract}
We prove the conjectural relation
between the Stokes matrix
for the quantum cohomology
of $X$
and an exceptional collection
generating $\DbX$
when $X$ is the Grassmannian
$\Gr(r,n)$.
The proof is based on the relation
between the quantum cohomology of the Grassmannian
and that of the projective space.
\end{abstract}

\section{Introduction}

Gromov-Witten invariants
of homogeneous spaces
contain enumerative information
such as the number of nodal rational curves
of a given degree
passing through a given set of points
in general position.
The theory of Frobenius manifold
allows a systematic treatment
of these invariants.
A Frobenius manifold is a complex manifold
whose tangent bundle has
a holomorphic bilinear form
and an associative commutative product
with certain compatibility conditions.
From these compatibility conditions,
it follows that there is
a function on the Frobenius manifold,
called the potential,
whose third derivatives give
the structure constants of the product.

Given a symplectic manifold $X$,
one can endow a Frobenius structure
on its total cohomology group $H^*(X;\bC)$.
In this case,
the holomorphic bilinear form
is given by the Poincar\'{e} pairing
and the potential is
the generating function
of the genus-zero Gromov-Witten invariants.
The product structure in this case
is called the quantum cohomology ring.
It is a deformation of the
cohomology ring
parametrized by $H^*(X;\bC)$
itself.

Given a Frobenius manifold,
one can construct the following isomonodromic family
of ordinary differential equations
on $\bP^1$:

\begin{equation} \label{eq:hbar_direction}
 \frac{\partial \Phi}{\partial \hbar}
  = (\frac{1}{\hbar}U+\frac{1}{\hbar^2} V) \Phi,
\end{equation}
\begin{equation} \label{eq:t_direction}
 \hbar \frac{\partial \Phi}{\partial t_\flata}
  = \frac{\partial}{\partial t_\flata} \circ \Phi, \quad
  \flata=0,\ldots,N-1.
\end{equation}

Here, $\Phi$ is the unknown function on $\bP^1$
times the Frobenius manifold
taking value in the tangent bundle of
the Frobenius manifold,
$\hbar$ is the coordinate on $\bP^1$,
$N$ is the dimension of the Frobenius manifold
and
$\{t_\flata\}_{\flata=0}^{N-1}$ is the {\em flat coordinate}
of the Frobenius manifold.
The circle denotes the product on the tangent bundle
and $U$, $V$ are certain operators
acting on sections of the tangent bundles.
See Dubrovin \cite{Dubrovin_G2DTFT}
for details.
Note that $z$ loc. cit. is
$1/\hbar$ in this paper.
(\ref{eq:hbar_direction}) is an
ordinary differential equation on $\bP^1$
with a regular singularity at infinity
and an irregular singularity at the origin,
and (\ref{eq:t_direction}) gives
its isomonodromic deformation.
If a point
on the Frobenius manifold
is semisimple,
i.e., if there are no nilpotent elements
in the product structure
on the tangent space
at this point,
one can define the {\em monodromy data}
of (\ref{eq:hbar_direction})
at this point,
consisting of the monodromy matrix at infinity,
the Stokes matrix at the origin
and the connection matrix between
infinity and the origin.
These data do not depend
on the choice of a semisimple point
because of the isomonodromicity.

The following conjecture,
originally due to Kontsevich,
developed by Zaslow \cite{Zaslow}, and
formulated into the following form
by Dubrovin \cite{Dubrovin_GATFM},
reveals a striking connection
between the Gromov-Witten invariants
and the derived category of coherent sheaves:

\begin{conjecture} \label{conj:stokes}

The quantum cohomology of a smooth projective variety $X$
is semisimple
if and only if
the bounded derived category $\DbX$
of coherent sheaves on $X$
is generated as a triangulated category
by an exceptional collection
$( \scE_i )_{i=1}^N$.
In such a case, the Stokes matrix $S$
for the quantum cohomology of $X$
is given by
\begin{equation} \label{eq:conjecture}
S_{ij} = \sum_k (-1)^k \dim \Ext^k(\scE_i, \scE_j).
\end{equation}
\end{conjecture}

An exceptional collection appearing above
is the following:

\begin{definition} \label{def:exceptional_collection}
\begin{enumerate}
\item 
 An object $\scE$ in a triangulated category
 is exceptional if
  $$
   \Ext^i(\scE,\scE)=\left\{\begin{array}{cl}
      \bC & \mbox{if $i=0$,}\\
       0  & \mbox{otherwise.}\\
   \end{array}\right.
  $$
\item
 An ordered set of objects $(\scE_i)_{i=1}^N$
 in a triangulated category
 is an exceptional collection
 if each $\scE_i$ is exceptional and
 $\Ext^k(\scE_i,\scE_j)=0$ for any $i>j$ and any $k$.

\end{enumerate}

\end{definition}

To our knowledge,
Conjecture \ref{conj:stokes} was previously known to hold
only for projective spaces
\cite{Dubrovin_PT2DTFT},
\cite{Guzzetti}.
The main result in this paper is:

\begin{theorem}
 Conjecture \ref{conj:stokes} holds
 for the Grassmannian $\Gr(r,n)$
 of $r$-dimensional subspaces
 in $\bC^n$.
\end{theorem}

The proof consists
of explicit computations
on both sides
of (\ref{eq:conjecture}).
The computation
on the left hand side
relies on the following two results:
The first is
a conjecture of Hori and Vafa
\cite{Hori-Vafa},
proved by Bertram, Ciocan-Fontanine and Kim
\cite{Bertram--Ciocan-Fontanine--Kim},
describing the solution
of (\ref{eq:hbar_direction})
for the Grassmannian $\Gr(r,n)$
in terms of that of the product of projective spaces
$(\bP^{n-1})^r$.
The second is the Stokes matrix
for the quantum cohomology of projective spaces
obtained by Dubrovin
\cite{Dubrovin_PT2DTFT}
for the projective plane
and by Guzzetti \cite{Guzzetti}
in any dimensions.
By combining these two results,
we can compute the Stokes matrix
for the quantum cohomology
of the Grassmannian.

On the right hand side,
we have an exceptional collection
generating
$D^b \coh(\Gr(r,n))$
by Kapranov \cite{Kapranov}.
It consists of equivariant vector bundles
on $\Gr(r,n)$
and $\Ext$-groups between them
can be computed
by the Borel-Weil theory.

Both of the above computations
can be carried out
for any $r$ and $n$,
and Conjecture \ref{conj:stokes}
reduces to the combinatorial identity
in Corollary \ref{cor:combinatorial}.

{\bf Acknowledgements}:
We thank A. N. Kirillov
for providing the proof
of the identity in Corollary \ref{cor:combinatorial}
and for allowing us to include it
in this paper.
We also thank H. Iritani,
T. Kawai, Y. Konishi, T. Maeno,
K. Saito and A. Takahashi
for valuable discussions and comments.
The author is supported by
JSPS Fellowships for Young Scientists
No.15-5561.

\section{Stokes matrix from the Hori-Vafa conjecture}

Let us begin with the discussion
of the Stokes matrix.
Fix a semisimple point
on a Frobenius manifold.
The differential equation
(\ref{eq:hbar_direction})
has a regular singularity at infinity
and an irregular singularity at the origin,
and the Stokes matrix is the monodromy data
for the irregular singularity at the origin,
defined as follows:
First, fix a formal fundamental solution
$\Phi_{\formal}$
of the form
\begin{equation} \label{eq:formal_solution}
 \Phi_{\formal}(\hbar)
  = \Psi R(\hbar) \exp[U / \hbar]
\end{equation}
where
$$
 U = \diag(u_1, \ldots, u_N),
$$
$\{u_\canonicala\}_{\canonicala=1}^N$
is the {\em canonical coordinate},
$\Psi$ is the coordinate transformation matrix
from the flat coordinate
to the normalized canonical coordinate
and
$R(\hbar) = (1+R_1 \hbar + R_2 \hbar^2+\cdots)$
is a formal series
satisfying
$$
 R^{t}(\hbar) R(-\hbar) = 1.
$$
Here, $\bullet^t$ denotes the transpose
of a matrix.
By \cite{Dubrovin_PT2DTFT}, Lemma 4.3.,
such $R(\hbar)$ exists uniquely.
Here we have taken the local trivialization
of the tangent bundle
given by the normalized canonical coordinate
and regarded $\Phi$
as an $n \times n$ matrix-valued function.

\begin{definition}
For $0 \leq \phi < \pi$,
a straight line
$l = \{\hbar \in \bCx \suchthat \arg(\hbar) = \phi, \phi - \pi \}$
passing through the origin
is called {\em admissible}
if the line through $u_k$ and $u_{k'}$
is not orthogonal to $l$
for any $k \neq k'$.
\end{definition}

Fix such a line, and choose
a small enough number $\epsilon > 0$
so that any line passing through the origin
with angle between $\phi-\epsilon$ and $\phi+\epsilon$
is admissible.
\begin{figure}[h]
 \centering
 \psfrag{Dl}{$D_\lt$}
 \psfrag{Dr}{$D_\rt$}
 \psfrag{D-}{$D_-$}
 \psfrag{phi}{$\phi$}
 \psfrag{l}{$l$}
 \includegraphics{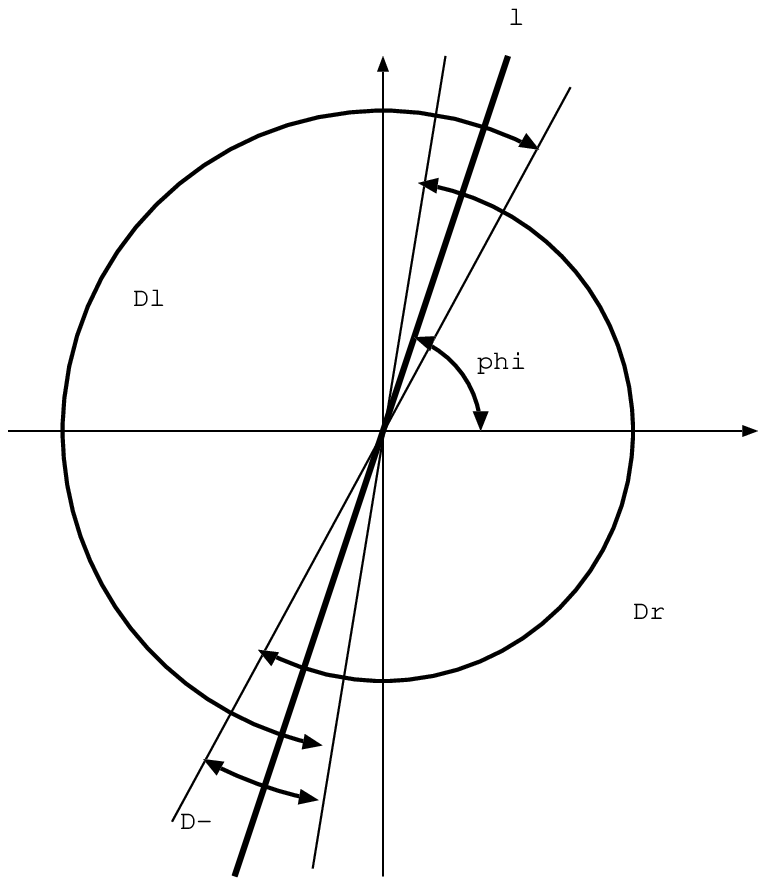}
\end{figure}
Define 
\begin{eqnarray} 
D_{\rt} & = & \{\hbar \in \bCx \suchthat
 \phi - \pi - \epsilon < \arg(\hbar) <\phi + \epsilon \},
 \nonumber \\
D_{\lt} & = & \{\hbar \in \bCx \suchthat
 \phi - \epsilon < \arg(\hbar) <\phi +\pi + \epsilon \}, \\
D_- & = & \{\hbar \in \bCx \suchthat
 \phi - \pi - \epsilon < \arg(\hbar) < \phi - \pi + \epsilon \}.
 \nonumber
\end{eqnarray}
Since the singularity at the origin
is irregular,
the formal solution $\Phi_\formal(\hbar)$
does not converge.
Nevertheless,
by \cite{Dubrovin_PT2DTFT}, Theorem 4.2.,
there exist unique solutions
$\Phi_\rt(\hbar)$ and $\Phi_\lt(\hbar)$,
defined on the angular domains $D_\rt$ and $D_\lt$
respectively,
which asymptote
to the same formal solution:
$$
\Phi_{\rt/\lt} \sim \Phi_{\small \mathrm{formal}}
\ \ \mbox{as $\hbar \rightarrow 0$ in $D_{\rt/\lt}$}.
$$
Since these two solutions
satisfy the same linear differential equation
on $D_-$,
there exists a matrix $S$
independent of $\hbar$
such that
\begin{equation*}
 \Phi_\rt(\hbar) = \Phi_\lt(\hbar) S,
 \qquad 
 \quad \hbar \in D_-.
\end{equation*}
This matrix $S$ is called the Stokes matrix.
Although locally on the Frobenius manifold
this Stokes matrix does not depend
on the choice of a semisimple point
by \cite{Dubrovin_PT2DTFT},
Isomonodromicity Theorem (second part),
it undergoes a discrete change
as we vary the point
on the Frobenius manifold
so that it crosses the point
where the line $l$
we have fixed at the beginning
is not admissible
any more.
This change in the Stokes matrix
is described
by an action of the braid group
$B_N$
(the number of strands
is the dimension of the Frobenius manifold).

In the case of the projective space $\bP^{n-1}$,
semisimplicity of the quantum cohomology
is well-known.
The solution to (\ref{eq:hbar_direction}),
(\ref{eq:t_direction}) has
an integral representation by Givental:
\begin{theorem}[Givental \cite{Givental_EGWI}] \label{th:Givental}
Let
$$
 W(x_1,\ldots,x_{n-1}) = x_1 + \cdots + x_{n-1}
  + \frac{e^t}{x_1 \cdots x_{n-1}}
$$
be a function on $(\bCx)^{n-1}$
depending on a parameter
$t \in \bC$
and choose a basis
$\{\Gamma_\canonicala\}_{\canonicala=1}^n$
of the space of flat sections
of the relative homology bundle
(the flat bundle on the $\hbar$-plane
whose fiber over $\hbar \in \bCx$ is
$H_{n-1}((\bCx)^{n-1}, \Re(W/\hbar) = -\infty)$).
Let $\{ p^\flata \}_{\flata=0}^{n-1}$
be the basis of $H^*(\bP^{n-1};\bZ)$
such that
$p^\flata \in H^{2 \flata}(\bP^{n-1};\bZ)$.
For $k=1,\ldots,n$, define
a cohomology-valued function $I_k$ by
\begin{equation} \label{eq:pn_integral_rep}
 I_k = \sum_{\flata = 0}^{n-1} p^\flata \int_{\Gamma_k}
  (\hbar \frac{d}{d t})^{\flata}
  \exp[{W(x_1,\ldots,x_{n-1})}/\hbar]
  \frac{d x_1 \cdots d x_{n-1}}{x_1 \cdots x_{n-1}}.
\end{equation}
Then $(I_k)_{k=1}^n$ gives a fundamental solution
to (\ref{eq:hbar_direction}), (\ref{eq:t_direction})
where $t$ is the coordinate of
$H^2(\bP^{n-1};\bC)$
and all the other flat coordinates
are set to zero.
\end{theorem}

Note that since the relative homology bundle
has a monodromy,
$\Gamma_k$'s
(and hence $I_k$'s)
cannot be defined globally.
The above integral representation
is related to the Stokes matrix
in the following way:
Fix $\phi$ and $\epsilon$ such that
any line passing through the origin
with angle between $\phi-\epsilon$ and $\phi+\epsilon$
is admissible.
There are $n$ critical points
and their critical values are
the canonical coordinate
$\{ u_\canonicala \}_{\canonicala=1}^n$.
Order these critical points
so that
$\Re[\exp(- \sqrt{-1} \phi) u_\canonicala]
 > \Re[\exp(- \sqrt{-1} \phi) u_\canonicalb]$
if $\canonicala < \canonicalb$.
Take the Lefschetz thimble
(the descending Morse cycle
for a suitable choice of a Riemannian metric
on $(\bCx)^{n-1}$)
for $\Re(W/\hbar)$
at $\hbar = \exp[\sqrt{-1}(\phi-\pi/2)]$
starting from the $\canonicala$-th critical point
of $W$
and extend it to a flat section
of the relative homology bundle
on $D_\rt$.
Let us call this section
$\Gamma_{\canonicala,\rt}$
and let $I_{\canonicala,\rt}$
be the integral
as in (\ref{eq:pn_integral_rep})
with $\Gamma_{\canonicala, \rt}$
as the integration cycle.
Now form the row vector
$(I_{\canonicala, \rt})_{\canonicala=1}^n$
and think of it
as an $n \times n$ matrix
by regarding an element
in the cohomology group
as a column vector
by the normalized canonical coordinate.
Then we can see that
$(I_{\canonicala, \rt})_{\canonicala=1}^n$
asymptotes on $D_{\rt}$
to the formal solution
of the form (\ref{eq:formal_solution})
as $\hbar \rightarrow 0$
by the saddle-point method.
In the same way,
starting from the Lefschetz thimble
at $\hbar = \exp[\sqrt{-1}(\phi+\pi/2)]$,
we obtain a solution
$(I_{\canonicala,\lt})_{\canonicala=1}^n$
defined on $D_{\lt}$
which asymptotes to the same formal solution
as $(I_{\canonicala,\rt})_{\canonicala=1}^n$.
Since the integrand is single-valued,
the monodromy of $I_\canonicala$'s comes
solely from the monodromy of the integration cycles
and the Stokes matrix is given by
$$
 \Gamma_{\canonicala, \rt}
  = \sum_{\canonicalb = 1}^n
     \Gamma_{\canonicalb, \lt} S_{\canonicalb \canonicala}.
$$

The Stokes matrix for the quantum cohomology
of $\bP^{n-1}$
has been computed
by Dubrovin \cite{Dubrovin_PT2DTFT} for $n \leq 3$
and by Guzzetti \cite{Guzzetti} for general $n$.
See also \cite{Tanabe_IHGAQCPS}.

\begin{theorem}[Dubrovin, Guzzetti] \label{th:Dubrovin_Guzzetti}
The Stokes matrix $S$
for the quantum cohomology
of the projective space $\bP^{n-1}$
is given by
$$
 S_{ij}= \combination{n-1+j-i}{j-i}
$$
up to the braid group action.
Here, $\combination{n}{r}$ is the binomial coefficient.
\end{theorem}

Since
$(\scO_{\bP^{n-1}}(i))_{i=0}^{n-1}$
is an exceptional collection
generating $D^b \coh \bP^{n-1}$
by Beilinson \cite{Beilinson}
and
$$
 \combination{n-1+j-i}{j-i}
  = \sum_k (-1)^k \dim \Ext^k(\scO_{\bP^{n-1}}(i), \scO_{\bP^{n-1}}(j)),
$$
the Conjecture \ref{conj:stokes}
holds for projective spaces.

Now let us move on
to the Grassmannian case.
Let $\Gr(r,n)$ be the Grassmannian
of $r$-dimensional subspaces
in $\bC^n$.
The semisimplicity of the quantum cohomology
in this case is also known.
The following theorem
is proved by
Bertram, Ciocan-Fontanine and Kim
(see the proof of Theorem 3.3
in \cite{Bertram--Ciocan-Fontanine--Kim}):

\begin{theorem}[Bertram--Ciocan-Fontanine--Kim]
\label{th:Hori-Vafa_conj}
For a choice of a basis
$\{ \phi_\flata \}_{\flata=0}^{N-1}$
of $H^*(\Gr(r,n);\bC)$
where
$N = \combination{n}{r}
   = \dim H^*(\Gr(r,n);\bC)$,
there exists a set
$\{ \varphi_\flata(x_{11}, \ldots, x_{r,n-1};t,\hbar) \}_{\flata=0}^{N-1}$
of functions of
$(x_{11}, \ldots, x_{r,n-1}) \in (\bC^\times)^{r(n-1)}$,
$t \in \bC$, and $\hbar \in \bCx$
such that
$$
\left( \sum_{\flata=0}^{N-1} \phi_\flata
 \int_{\Gamma_{k_1} \times \cdots \Gamma_{k_r}} e^{W/\hbar}
  \varphi_\flata(x_{11}, \ldots, x_{r,n-1};t,\hbar)
  \prod_{j=1}^r
  \frac{d x_{j1} \cdots d x_{j, n-1}}{x_{j1} \cdots x_{j, n-1}}
\right)_{1 \leq k_1 < k_2 < \cdots < k_r \leq n}
$$
forms a fundamental solution
to (\ref{eq:hbar_direction}), (\ref{eq:t_direction})
where $t$ is the coordinate of
$H^2(\Gr(r,n);\bC)$
and all the other flat coordinates
are set to zero.
Here,
$$
 W(x_{11}, \ldots, x_{r,n-1})
  = \sum_{j=1}^r
     \left( 
      x_{j1} + \cdots + x_{j,n-1}
      + \frac{e^t}{x_{j1} \cdots x_{j,n-1}}
     \right)
$$
and $\{\Gamma_k\}_{k=1}^n$
is the basis of the flat sections
of the relative homology bundle
as in Theorem \ref{th:Givental}.
\end{theorem}

By construction,
$\varphi_\flata(x_{11}, \ldots, x_{r,n-1};t,\hbar)$
is anti-symmetric
with respect to the exchange of
$(x_{i1}, \cdots, x_{i,n-1})$
and $(x_{j1}, \cdots, x_{j,n-1})$
for any $1 \leq i < j \leq r$.
Therefore,
if we define
$H^*(\Gr(r,n);\bC)$-valued functions
$I_K(t,\hbar)$
for $K=(k_1, \ldots, k_r)$,
$1 \leq k_i \leq n$,
$i = 1, \ldots, r$
by
$$
 I_K = \sum_{\flata=0}^{N-1} \phi_\flata
 \int_{\Gamma_{k_1} \times \cdots \Gamma_{k_r}} e^{W/\hbar}
  \varphi_\flata(x_{11}, \ldots, x_{r,n-1};t,\hbar)
  \prod_{j=1}^r
  \frac{d x_{j1} \cdots d x_{j, n-1}}{x_{j1} \cdots x_{j, n-1}},
$$
then $I_K$ is totally anti-symmetric
in $k_1, \ldots, k_r$.
Hence it follows that
if we put
$$
 \Gamma_K = \frac{1}{r!} \sum_{\sigma \in \frakS_r}
  \sgn \sigma \; \Gamma_{k_{\sigma(1)}} \times \cdots
   \times \Gamma_{k_{\sigma(r)}}
$$
where $\frakS_r$ is the symmetric group
of degree $r$
and $\sgn \sigma$ is the signature of $\sigma$,
then we have
$$
 I_K = \sum_{\flata=0}^{N-1} \phi_\flata
 \int_{\Gamma_K} e^{W/\hbar}
  \varphi_\flata(x_{11}, \ldots, x_{r,n-1};t,\hbar)
  \prod_{j=1}^r
  \frac{d x_{j1} \cdots d x_{j, n-1}}{x_{j1} \cdots x_{j, n-1}}.
$$
We can use the above result
to compute the Stokes matrix
for the quantum cohomology of $\Gr(r,n)$
from that of $\bP^{n-1}$ as follows:
By Theorem \ref{th:Dubrovin_Guzzetti},
there exists a choice
$\{\Gamma_{\canonicala,\rt}\}_{\canonicala=1}^n$ and 
$\{\Gamma_{\canonicala,\lt}\}_{\canonicala=1}^n$
of bases
of flat sections of the relative homology bundle
on $D_\rt$ and $D_\lt$ respectively
such that
$$
 \Gamma_{\canonicala, \rt}
  = \sum_{\canonicalb = 1}^n
     \Gamma_{\canonicalb, \lt} S_{\canonicalb \canonicala}
$$
on $D_-$ for $S_{ij}= \combination{n-1+j-i}{j-i}$.
Then the monodromy for
$\Gamma_K$ is given by
$$
 \Gamma_{K,\rt}
  = \sum_{1 \leq l_1 < l_2 < \cdots < l_r \leq n}
       I_{L,\lt} \; S_{L, K}.
$$
where $K=(k_1,\ldots,k_r)$,
$L=(l_1,\ldots,l_r)$, and
\begin{eqnarray}
S_{L, K}
  & = & \det (S_{l_j, k_i})_{1 \leq i,j \leq r} \nonumber \\
  & = & \det(\combination{n+l_i-k_j-1}{l_i-k_j})_{1 \leq i,j \leq r}.
    \label{eq:stokes}   
\end{eqnarray}

\section{Derived category of coherent sheaves}

In this section, we use the presentation
$$
 \Gr(r,n) = \GL_n(\bC) / P
$$
of the Grassmannian as a homogeneous space,
where
$$ 
 P = \left\{\begin{pmatrix}
	     A & B \\ 0 & D
	    \end{pmatrix} \suchthat A \in \GL_r(\bC),
  B \in M_{r, n-r}(\bC), D \in \GL_{n-r}(\bC) \right\}
$$
is a parabolic subgroup of $GL_n(\bC)$.
A representation of $\GL_r(\bC)$ gives
a representation of $P$
through the projection
$P \ni \begin{pmatrix}
 A & B \\ 0 & D \end{pmatrix}
   \mapsto A \in \GL_r(\bC)$,
hence a $\GL_n(\bC)$-equivariant bundles
on $\Gr(r,n)$
associated to the principal $P$-bundle
$\GL_n(\bC) \rightarrow \Gr(r,n)$.
Let $\scE_\rho$ denote
the equivariant bundle on $\Gr(r,n)$
corresponding to a representation
$\rho$ of $\GL_r(\bC)$
in this way.

Let
$$
 \Lambda = \{(\lambda_1,\ldots,\lambda_r) \in \bZ^r
  \suchthat n-r \geq \lambda_1 \geq
   \cdots \geq \lambda_r \geq 0 \}
$$
be a set of weights of $\GL_r(\bC)$.
Given a weight $\lambda$,
let $\rho_\lambda$ denote
the irreducible representation
of $\GL_r(\bC)$
with highest weight $\lambda$.
We abbreviate $\scE_{\rho_\lambda}$
as $\scE_\lambda$.

%

\begin{theorem}[Kapranov \cite{Kapranov}]
 $\{\scE_{\lambda}\}_{\lambda \in \Lambda}$
 is an exceptional collection
 generating $D^b \coh(\Gr(r,n))$.
\end{theorem}

Kapranov also proved that
$\Ext^k(\scE_{\lambda},\scE_{\mu})=0$
for any $\lambda, \mu \in \Lambda$
and any $k \neq 0$.
$\Hom(\scE_{\lambda},\scE_{\mu})$
is calculated as follows:
Decompose the tensor product
$\rho_\lambda^\vee \otimes \rho_\mu$
of the dual representation of $\rho_\lambda$
and $\rho_\mu$
into the direct sum of
irreducible representations
$$
\rho_\lambda^\vee \otimes \rho_\mu
 = \bigoplus_{\nu} \rho_\nu^{\oplus \widetilde{N}_{\lambda \mu}^\nu}.
$$
Here, $\widetilde{N}_{\lambda \mu}^\nu$ is the multiplicity
of $\rho_\nu$ in $\rho_\lambda^\vee \otimes \rho_\mu$
and $\nu$ runs over all weights
of $\GL_r(\bC)$.
Define
$$
N_{\lambda \mu}^\nu = \left\{\begin{array}{cl}
      \widetilde{N}_{\lambda \mu}^\nu & \mbox{if $\nu_r \geq 0$,}\\
       0  & \mbox{otherwise.}\\
   \end{array}\right.
$$
For a weight $\lambda \in \bZ^r$ of $\GL_r(\bC)$,
let $R_\lambda$ be the irreducible representation
of $\GL_n(\bC)$ with highest weight
$(\lambda_1,\ldots,\lambda_r,0,\ldots,0) \in \bZ^n$.
Then
\begin{eqnarray}
 \Hom(\scE_{\lambda},\scE_{\mu})
 & = & H^0(\scE_{\lambda}^\vee \otimes \scE_{\mu}) \nonumber \\
 & = & H^0(\scE_{\rho_\lambda^\vee \otimes \rho_\mu}) \nonumber \\
 & = & \bigoplus_{\nu} H^0(\scE_{\nu})^
         {\oplus \widetilde{N}_{\lambda \mu}^\nu} \nonumber \\
 & = & \bigoplus_{\nu} R_{\nu}^{\oplus N_{\lambda \mu}^\nu}, \label{eq:DbX}
\end{eqnarray}
where the last equality follows
from the Borel-Weil theory.

\section{A combinatorial identity} \label{section:combinatorics}


The content of this section
is due to A. N. Kirillov.
Fix two integers $r$, $n$ such that $r < n$.
Let $A=\{(k_1,\ldots,k_r) \in \bZ^r \suchthat
1 \leq k_1 < \cdots < k_r \leq n \}$.
$A$ and $\Lambda$ defined in the previous section
are bijective by the correspondence
$$
\Lambda \ni (\lambda_i)_{i=1}^r
 \mapsto (\lambda_{r-i+1}+i)_{i=1}^r \in A.
$$
For $n$ variables $x=(x_1,\ldots,x_n)$,
let $s_\lambda(x) = \det(h_{\lambda_i-i+j}(x))_{1 \leq i,j \leq n}$
be the Shur function,
where $h_i(x)$ is the complete symmetric function
(the sum of all monomials of degree $i$).
For generalities on
symmetric functions,
see, e.g.,
\cite{Macdonald}.
Define integers $c_{\mu \nu}^\lambda$'s by
$$
s_\mu(x) s_\nu(x) = \sum_\lambda c_{\mu \nu}^\lambda s_\lambda(x)
$$
and the skew Shur function $s_{\lambda/\mu}(x)$ by
$$
s_{\lambda/\mu}(x) = \sum_{\nu} c_{\mu \nu}^\lambda s_\nu(x).
$$
Then
$$
s_{\lambda/\mu}(x) = \det(h_{\lambda_i-\mu_j-i+j}(x))_{1 \leq i,j \leq n}.
$$
\begin{lemma}
Let $\mu$, $\nu$ and $\lambda$ be partitions
such that $\mu_1 \leq \nu_r$.
Define $\mu^c = (\mu_1-\mu_r, \mu_1-\mu_{r-1},\ldots,\mu_1-\mu_2,0)$
and $\widetilde{\nu} = (\nu_1-\mu_1, \nu_2-\mu_1, \ldots, \nu_r-\mu_1)$.
Then
$$
c_{\lambda \mu^c}^\nu = c_{\mu \widetilde{\nu}}^\lambda.
$$
\end{lemma}
\proof
\begin{eqnarray*}
 c_{\lambda \mu^c}^\nu
  &=& \dim \Hom_{\GL_r(\bC)}(\rho_\lambda \otimes \rho_{\mu^c}, \rho_\nu) \\
  &=& \dim \Hom_{\GL_r(\bC)}(\rho_0, \rho_\lambda^\vee
        \otimes \rho_{\mu^c}^\vee \otimes \rho_\nu) \\
  &=& \dim \Hom_{\GL_r(\bC)}(\rho_0, \rho_\lambda^\vee
        \otimes (\rho_\mu^\vee \otimes \det^{\otimes \mu_1})^\vee
        \otimes \rho_\nu) \\
  &=& \dim \Hom_{\GL_r(\bC)}(\rho_0, \rho_\lambda
        \otimes \rho_\mu^\vee \otimes \det^{\otimes \mu_1}
        \otimes \rho_\nu^\vee) \\
  &=& \dim \Hom_{\GL_r(\bC)}(\rho_0, \rho_\lambda
        \otimes \rho_\mu^\vee
        \otimes \rho_{\widetilde{\nu}}^\vee) \\
  &=& \dim \Hom_{\GL_r(\bC)}(\rho_\mu \otimes \rho_{\widetilde{\nu}},
        \rho_\lambda),
\end{eqnarray*}
where $\rho_0$ is the trivial representation
and $\det$ is the determinant representation
(the irreducible representation
with highest weight $(1,\ldots,1) \in \bZ^r$).
\qed

\begin{theorem} \label{th:combinatorial}
 $\displaystyle
   s_{\lambda/\mu}(x)= \sum_\nu N_{\mu \lambda}^\nu s_\nu(x).$
\end{theorem}
\proof
\begin{eqnarray*}
\sum_\nu N_{\mu \lambda}^\nu s_\nu(x)
 &=& \sum_\nu c_{\lambda \mu^c}^\nu s_{\widetilde{\nu}}(x) \\
 &=& \sum_{\widetilde{\nu}} c_{\mu \widetilde{\nu}}^\lambda
                      s_{\widetilde{\nu}}(x) \\
 &=& s_{\lambda/\mu}(x).
\end{eqnarray*}
\qed

By substituting $x_1=\cdots=x_n=1$ in Theorem \ref{th:combinatorial}
and using $h_r(1,\ldots,1)=\combination{n+r-1}{r}$,
we obtain the following:

\begin{corollary} \label{cor:combinatorial}
For $\lambda, \mu \in \Lambda$, let
$k=(\lambda_{r-i+1}+i)_{i=1}^r$,
$l=(\mu_{r-i+1}+i)_{i=1}^r$.
Then
$$
 \det(\combination{n+l_i-k_j-1}{l_i-k_j})_{1 \leq i,j \leq r}
  = \sum_\nu N_{\lambda \mu}^\nu \dim R_\nu.
$$
\end{corollary}

The left hand side is the component
of the Stokes matrix from (\ref{eq:stokes})
and the right hand side is the Euler number
in the derived category of coherent sheaves
form (\ref{eq:DbX}).
This proves Conjecture \ref{conj:stokes}
in the case of Grassmannians.

\bibliographystyle{plain}
\bibliography{bibs}

\newcommand{\noop}[1]{}\def\cprime{$'$} \def\cprime{$'$}
\begin{thebibliography}{10}

\bibitem{Beilinson}
A.~A. Be{\u\i}linson.
\newblock Coherent sheaves on {${\bf P}\sp{n}$} and problems in linear algebra.
\newblock {\em Funktsional. Anal. i Prilozhen.}, 12(3):68--69, 1978.

\bibitem{Bertram--Ciocan-Fontanine--Kim}
A.~Bertram, I.~Ciocan-Fontanine, and B.~Kim.
\newblock Two proofs of a conjecture of {H}ori and {V}afa.
\newblock math.AG/0304403.

\bibitem{Dubrovin_G2DTFT}
Boris Dubrovin.
\newblock Geometry of {$2$}{D} topological field theories.
\newblock In {\em Integrable systems and quantum groups (Montecatini Terme,
  1993)}, volume 1620 of {\em Lecture Notes in Math.}, pages 120--348.
  Springer, Berlin, 1996.

\bibitem{Dubrovin_GATFM}
Boris Dubrovin.
\newblock Geometry and analytic theory of {F}robenius manifolds.
\newblock In {\em Proceedings of the International Congress of Mathematicians,
  Vol. II (Berlin, 1998)}, number Extra Vol. II, pages 315--326 (electronic),
  1998.

\bibitem{Dubrovin_PT2DTFT}
Boris Dubrovin.
\newblock Painlev\'e transcendents in two-dimensional topological field theory.
\newblock In {\em The Painlev\'e property}, CRM Ser. Math. Phys., pages
  287--412. Springer, New York, 1999.

\bibitem{Givental_EGWI}
Alexander~B. Givental.
\newblock Equivariant {G}romov-{W}itten invariants.
\newblock {\em Internat. Math. Res. Notices}, (13):613--663, 1996.

\bibitem{Guzzetti}
Davide Guzzetti.
\newblock Stokes matrices and monodromy of the quantum cohomology of projective
  spaces.
\newblock {\em Comm. Math. Phys.}, 207(2):341--383, 1999.

\bibitem{Hori-Vafa}
K.~Hori and C.~Vafa.
\newblock Mirror symmetry.
\newblock hep-th/0002222.

\bibitem{Kapranov}
M.~M. Kapranov.
\newblock On the derived categories of coherent sheaves on some homogeneous
  spaces.
\newblock {\em Invent. Math.}, 92(3):479--508, 1988.

\bibitem{Macdonald}
I.~G. Macdonald.
\newblock {\em Symmetric functions and {H}all polynomials}.
\newblock Oxford Mathematical Monographs. The Clarendon Press Oxford University
  Press, New York, second edition, 1995.
\newblock With contributions by A. Zelevinsky, Oxford Science Publications.

\bibitem{Tanabe_IHGAQCPS}
Susumu Tanab{\'e}.
\newblock Invariant of the hypergeometric group associated to the quantum
  cohomology of the projective space.
\newblock {\em Bull. Sci. Math.}, 128(10):811--827, 2004.

\bibitem{Zaslow}
Eric Zaslow.
\newblock Solitons and helices: the search for a math-physics bridge.
\newblock {\em Comm. Math. Phys.}, 175(2):337--375, 1996.

\end{thebibliography}

\end{document}